# Dynamics & Stochastics

## Festschrift in honor of M. S. Keane

Dee Denteneer, Frank den Hollander, Evgeny Verbitskiy, Editors

# Contents

# Contributors to this volume

Baggen, S., *Philips Research Laboratories Eindhoven*
Balakirsky, V., *Eindhoven University of Technology*
Beffara, V., *CNRS–ENS-Lyon*
Bosma, W., *Radboud University Nijmegen*
Broer, H. W., *Groningen University*
Bulinski, A., *Moscow State University*

Coven, E. M., *Wesleyan University*

Dajani, K., *Utrecht University*
Dekking, F. M., *Thomas Stieljes Institute for Mathematics and Delft University of Technology*
den Hollander, F., *Leiden University & EURANDOM*
Denteneer, D., *Philips Research Laboratories Eindhoven*

Egner, S., *Philips Research Laboratories Eindhoven*

Grimmett, G., *University of Cambridge*
Guivarc'h, Y., *Université de Rennes 1*

Häggström, O., *Chalmers University of Technology*
Hanßmann, H., *RWTH Aachen, Universiteit Utrecht*
Hollmann, H., *Philips Research Laboratories Eindhoven*

Kümmerer, B., *Technische Universität Darmstadt*
Kahn, J., *Rutgers University*
Kamae, T., *Matsuyama University*
Keane, M., *Wesleyan University*
Kraaikamp, C., *University of Technology Delft*
Krieger, W., *University of Heidelberg*

Liardet, P., *Université de Provence*

Maass, A., *Universidad de Chile*
Maassen, H., *Radboud University Nijmegen*
Martínez, S., *Universidad de Chile*
Meester, R., *Vrije Universiteit*
Merkl, F., *University of Munich*
Misiewicz, J. K., *University of Zielona Góra*
Mori, M., *Nihon University*

Naudot, V., *Groningen University*
Netočný, K., *Institute of Physics*

Petersen, K., *University of North Carolina*
Pivato, M., *Trent University*

Queffélec, M., *Université Lille1*

Rolles, S. W. W., *University of Bielefeld*
Roy, R., *Indian Statistical Institute*

Schaap, H. G., *University of Groningen*
Schmidt, K., *University of Vienna and Erwin Schrödinger Institute*
Schweiger, F., *University of Salzburg*
Serafin, J., *Wrocław University of Technology*
Shashkin, A., *Moscow State University*
Sidoravicius, V., *IMPA*
Spohn, H., *TU-München*
Steif, J. E., *Chalmers University of Technology*

Tolhuizen, L., *Philips Research Laboratories Eindhoven*

van den Berg, J., *CWI and VUA*
van der Poorten, A. J., *Centre for Number Theory Research, Sydney*
van Enter, A. C. D., *University of Groningen*
van Lieshout, M. N. M., *Centre for Mathematics and Computer Science, Amsterdam*
Vares, M. E., *CBPF*
Verbitskiy, E., *Philips Research Laboratories Eindhoven*

Ward, T., *University of East Anglia*

Yassawi, R., *Trent University*

# Preface

The present volume is a *Festschrift* for Mike Keane, on the occasion of his 65th birthday on January 2, 2005. It contains 29 contributions by Mike's closest colleagues and friends, covering a broad range of topics in *Dynamics and Stochastics*.

To celebrate Mike's scientific achievements, a conference entitled "Dynamical Systems, Probability Theory and Statistical Mechanics" was organized in Eindhoven, The Netherlands, during the week of January 3–7, 2005. This conference was hosted jointly by EURANDOM and by Philips Research. It drew over 80 participants from 5 continents, which is a sign of the warm affection and high esteem for Mike felt by the international mathematics community.

Mike is one of the founders of EURANDOM, and has seen this institute come to flourish since it opened its doors in 1998, bringing much extra buzz and liveliness to Dutch stochastics. Mike has also been a scientific consultant to Philips Research for some 20 years, and since 1998 spends one month per year with Philips Research Eindhoven in that capacity. It is therefore particularly nice that the conference could take place under the umbrella of these two institutions so close to Mike's heart.

Most people retire at 65. Not Mike, who continues to be enormously energetic and full of plans. After a highly active career, taking him to professorships in a number of countries (Germany 1962–1968, USA 1968–1970, France 1970–1980, The Netherlands 1981–2002), Mike is now with Wesleyan University in Middletown, Connecticut, USA, from where he continues his relentless search for fundamental questions, elegant solutions and powerful applications. All this takes place in an atmosphere of warm hospitality through the help of Mike's spouse Mieke, who has welcomed more guests at her home than the American ambassador.

A special event took place in the late afternoon of January 5, 2005, in the main auditorium on the Philips High Tech Campus. On behalf of Her Majesty the Queen, Mr. A. B. Sakkers, the mayor of Eindhoven, presented Mike with the decoration in the highest civilian order: "Knight in the Order of the Dutch Lion." This decoration was bestowed upon Mike for his outstanding scientific contributions, his stimulating international role, as well as his service to industry.

A number of organizations generously supported the conference:

- NWO (Netherlands Organisation for Scientific Research).
- KNAW (Royal Netherlands Academy of Arts and Sciences).
- MRI (Mathematical Research Institute) and TSI (Thomas Stieltjes Institute), the two main collaborative research schools for mathematics in The Netherlands.
- University of Amsterdam: Korteweg de Vries Institute.
- Eindhoven University of Technology: Department of Mathematics and Computer Science.
- Wesleyan University: Department of Mathematics.
- ESF (European Science Foundation), through its scientific program "Random Dynamics of Spatially Extended Systems."

We are grateful to these organizations for their contribution. We are also grateful to the program committee (Rob van den Berg, Michel Dekking, Hans van Duijn, Chris Klaassen, Hans Maassen and Ronald Meester) for their help with putting together the conference program.

We wish Mike and Mieke many healthy and active years to come, and trust that the reader will find in this Festschrift much that is to her/his interest and liking.

Eindhoven, December 2005

Dee Denteneer
Frank den Hollander
Evgeny Verbitskiy